\RequirePackage{amsmath}
\documentclass{aip-cp}

\usepackage[numbers]{natbib}
\usepackage{rotating}
\usepackage{graphicx}
\usepackage{verbatim}
\usepackage[normalem]{ulem} 
\newcommand{\remove}[1]{\blue{\sout{#1}}}

\newcommand{\N}{{\mathbb N}}
\newcommand{\R}{{\mathbb R}}
\newcommand{\C}{{\mathbb C}}

\usepackage{xcolor}

\newcommand{\blue}{\textcolor{blue}}

\newtheorem{theorem}{Theorem}
\newtheorem{lemma}[theorem]{Lemma}

\newtheorem{proposition}[theorem]{Proposition}

\newtheorem{definition}{Definition}

\newcommand{\Cnn}{\C^{n\times n}}

\newcommand{\mHC}{{m_{HC}}}

\newcommand{\mLH}{\mathbf{L}_H}
\newcommand{\mLS}{\mathbf{L}_S}

\newcommand{\mI}{\mathbf{I}}

\newcommand{\mL}{\mathbf{L}}
\newcommand{\mM}{\mathbf{M}}
\newcommand{\mP}{\mathbf{P}}

\newcommand{\dd}[1][x]{\,{\rm d} #1}
\newcommand{\dt}{\,\dd[t]}

\newcommand{\bigO}{\mathcal{O}}

\newcommand{\rank}{\rm rank\,}
\DeclareMathOperator{\sgn}{sgn}
\DeclareMathOperator{\spn}{span}

\begin{document}

\title{Hypocoercivity for Linear ODEs and Strong Stability for Runge--Kutta Methods}

\author[aff1]{Franz Achleitner}
\eaddress{franz.achleitner@tuwien.ac.at}
\author[aff1]{Anton Arnold\corref{cor1}}
\eaddress[url]{https://www.asc.tuwien.ac.at/arnold/}
\author[aff1]{Ansgar J\"ungel}
\eaddress{ansgar.juengel@tuwien.ac.at}

\affil[aff1]{Inst.\ f.\ Analysis u.\ Scientific Computing, Technische Universit\"at Wien,  Wiedner Hauptstr. 8, A-1040 Wien, Austria.}
\corresp[cor1]{Corresponding author: anton.arnold@tuwien.ac.at}

\maketitle

\begin{abstract}
In this note, we connect two different topics from linear algebra and numerical analysis: hypocoercivity of semi-dissipative matrices and strong stability for explicit Runge--Kutta schemes. Linear autonomous ODE systems with a non-coercive matrix are called hypocoercive if they still exhibit uniform exponential decay towards the steady state. Strong stability is a property of time-integration schemes for ODEs that preserve the temporal monotonicity of the discrete solutions. It is proved that explicit Runge--Kutta schemes are strongly stable with respect to semi-dissipative, asymptotically stable matrices if the hypocoercivity index is sufficiently small compared to the order of the scheme. Otherwise, the Runge--Kutta schemes are in general not strongly stable. As a corollary, explicit Runge--Kutta schemes of order $p\in 4\N$ with $s=p$ stages turn out to be \emph{not} strongly stable. This result was proved in \cite{AAJ23}, filling a gap left open in \cite{SunShu19}. Here, we 
present an alternative, direct proof.
\end{abstract}


\section{Hypocoercive ODEs}

For linear autonomous ordinary differential equations (ODEs),
\begin{equation}\label{ODE}
  \frac{\dd[u]}{\dt} =\mL u,\quad t>0, \quad u(0)=u^0\in\C^n,
\end{equation}
with Lyapunov stable matrices $\mL\in\C^{n\times n}$ (i.e., all eigenvalues of $\mL$ have nonpositive real part and the purely imaginary eigenvalues are non-defective), we are concerned with characterizing their short-time decay behavior. To this end, we review first hypocoercivity properties of such systems \cite{AAC18, Vi09}:

\begin{definition}\label{def:HC}
\begin{enumerate}
\item[(a)] The matrix $\mL\in\C^n$ is called \emph{dissipative} (resp.\ \emph{semi-dissipative}) if its Hermitian part, $\mLH:=(\mL+\mL^*)/2$ is negative definite (resp.\ negative semi-definite).
\item[(b)] $-\mL$ is called \emph{hypocoercive} (or \emph{positive stable}) if there are constants $\lambda>0$ and $c\ge1$ such that the matrix exponential satisfies
$$
  \| e^{\mL t}\|_2 \le c e^{-\lambda t},\quad t\ge0.
$$
\item[(c)] Let $\mL\in\C^n$ be semi-dissipative. Its \emph{hypocoercivity index (HC-index)} $\mHC$ is defined as the smallest integer $m\in\N_0$ (if it exists) such that 
$$
  T_m:=\sum_{j=0}^m \mLS^j \mLH (\mLS^*)^j <0, 
$$
where $\mLS:=(\mL-\mL^*)/2$ denotes the skew-Hermitian part of $\mL$.
\end{enumerate}
\end{definition}
The HC-index of $\mL$ characterizes the structural complexity of the interplay between $\mLH$ and $\mLS$, and it is bounded by 
$$
  0\le\frac{n-\rank \mLH}{\rank \mLH} \le \mHC(\mL) \le n-{\rank \mLH}
  \le n - 1;
$$
see \cite{AAC18, AAC23, AAM23}.
Hence the matrix $\mL$ is dissipative if and only if $\mHC(\mL)=0$. For later use, we consider the example
\begin{equation}\label{L}
 \mL
=\left( \begin{array}{rrrrr}
	0  & -1 & & & \\
	1 & \ddots & \ddots & & \\
	   & \ddots & \ddots & -1 & \\
		 & & 1 & 0 & -1 \\
		 & & & 1 & -1
	\end{array}\right) \in \R^{N\times N}.
\end{equation}
Using the Kalman rank condition \cite[Proposition 1]{AAC18}, one easily verifies that $\mL$ has (maximal) HC-index $N-1$. Moreover, the matrix $\mL$ is asymptotically stable, i.e., it has only eigenvalues with negative real part.

We recall that the HC-index gives a precise characterization of the short-time behavior of the solutions to \eqref{ODE}:

\begin{proposition}[\cite{AAC22}] 
 \label{prop:ODE-short}
Let the matrix $\mL\in\Cnn$ be semi-dissipative.
Then~$\mL$ is asymptotically stable (with HC-index $\mHC\in\N_0$) if and only if
\begin{equation}\label{short-t-decay}
  \|e^{t\mL}\|_2 = 1-ct^a+\mathcal{O}(t^{a+1})\quad\mbox{ for } t\in[0,\varepsilon), 
\end{equation}
for some $a,c,\varepsilon>0$. In this case, necessarily $a=2m_{HC}+1$.
\end{proposition}

The sharp multiplicative factor~$c$ in~\eqref{short-t-decay} has been determined explicitly in~\cite[Theorem 2.7(b)]{AAC22}.


\section{Strong stability for Runge--Kutta methods}

It is known that a matrix~$\mL\in\Cnn$ is \emph{Lyapunov stable} if and only if there exists a positive definite Hermitian matrix $\mP$ such that  
\begin{equation} \label{ineq:Lyapunov}
 \mL^* \mP +\mP\mL\leq 0.
\end{equation}  
In this case, the solution $u(t)$ of \eqref{ODE} is nonincreasing in the norm $\|\cdot\|_\mP:=\sqrt{\langle\cdot, \mP\cdot\rangle}$~:
\begin{equation} \label{ineq.monotonicity}
  \frac{\dd[]}{\dt}\|u(t)\|_\mP^2
= \langle u,(\mL^*\mP+\mP\mL)u\rangle 
\le 0.
\end{equation}
It is often desirable that a numerical scheme for \eqref{ODE} reproduces this decay behavior on the discrete level. 

Following \cite{AAJ23}, we consider here explicit Runge--Kutta schemes, where $u^k$ is an approximation for $u(k\tau),\,k\in\N_0$, and $\tau$ is the uniform time step:
\begin{equation}\label{1.RK}
  u^k 
= u^{k-1} + \tau\sum_{i=1}^s b_i K_i^k, \quad
  K_i^k
= \mL\bigg(u^{k-1} + \tau\sum_{j=1}^{i-1}a_{ij}K_j^k\bigg), \quad
	i=1,\ldots,s,
\end{equation}
where $b_i\in\C$ are the weights, $a_{ij}\in\C$ are the coefficients of the Runge--Kutta matrix, and $s\in\N$ is the number of stages. This scheme can be rewritten in compact form as $u^k=R(\tau\mL)u^{k-1}$, using its \emph{stability function} $R(z)$. For the scheme \eqref{1.RK} to have order $p$, one needs at least $s\ge p$ stages. In this case, its stability function takes the form
\begin{equation}\label{1.R}
 R(z) 
=\sum_{j=0}^p\frac{z^j}{j!} + \sum_{j=p+1}^s c_j\frac{z^j}{j!}, 
\quad z\in\C,\quad c_{p+1}\neq 1 ,
\end{equation}
with some constants $c_{p+1},\ldots,c_s\in\C$.

For the discrete analog of the monotonicity estimate  \eqref{ineq.monotonicity}, we shall use the following notions:

\begin{definition}\label{def.stab}
\begin{itemize}
\item[(a)] 
The Runge--Kutta scheme \eqref{1.RK} is {\em strongly stable} if for all matrix dimensions $n\in\N$, for all Lyapunov stable matrices $\mL\in\Cnn$, and for all Hermitian matrices $\mP>0$ such that~\eqref{ineq:Lyapunov} holds, the numerical solution to \eqref{ODE} satisfies $\|u^1\|_\mP\le\|u^0\|_\mP$ for all initial data $u^0\in\C^n$ and sufficiently small time steps.
\item[(b)] 
The Runge--Kutta scheme \eqref{1.RK} is {\em strongly stable w.r.t.\ a subset $\mathcal L_0$} of Lyapunov stable matrices (of any dimension $n$), if the condition from (a) holds for all $\mL\in\mathcal L_0$.
\end{itemize}
\end{definition}

Strong stability of explicit Runge--Kutta schemes was studied in, e.g., \cite[\S4]{Tad02} and \cite{SunShu19}. We shall focus on explicit Runge--Kutta schemes with $s=p$ stages. Their stability function is given by the first sum in \eqref{1.R}, and their stability behavior was analyzed in \cite{SunShu19}:
$$
  \mbox{They are } \quad 
  \left\{
  \begin{array}{ll}
  \mbox{strongly stable} & \quad \mbox{if }p\in4\N_0+3\ ,\\
  \mbox{not strongly stable} & \quad \mbox{if } p\in4\N_0+1 \mbox{ or } p\in4\N_0+2\ , \\
  \end{array}
  \right.
$$
but the following instability result for the case $p\in 4\N$ was only found and proved recently in \cite{AAJ23}:

\begin{theorem}[\cite{AAJ23}]	\label{th:p4N}
Explicit Runge--Kutta schemes of order $p\in4\N$ with $s=p$ stages are \emph{not} strongly stable.
\end{theorem}

While the proof of this result in \cite{AAJ23} was based on the quite technical result in Proposition \ref{prop:ODE-short}, we shall give here an independent direct proof. In order to motivate our subsequent proof, we first cite another result from \cite{AAJ23}. To this end, we define the following subset of asymptotically stable (and thus Lyapunov stable) matrices:
$$
  \mathcal L_{AS}^m:=\{\mL \mbox{ is semi-dissipative and asymptotically stable} : \mHC(\mL)\le m\},\quad m\in\N_0.
$$

\begin{proposition}[\cite{AAJ23}]\label{prop:LASm}
All explicit Runge--Kutta schemes of order $p\in \N$ (with $s\ge p$ stages) are strongly stable w.r.t.~$\mathcal L_{AS}^{m}$, if $m\in\N_0$ satisfies $2m+1\leq p$. 
\end{proposition}

\medskip
\noindent\emph{Proof (of Theorem \ref{th:p4N}).}
For each fixed $p\in4\N$, we shall analyze an asymptotically stable matrix $\mL$ as a counterexample. Due to Proposition \ref{prop:LASm}, the HC-index of such $\mL$ must be ``large enough''. More precisely, we choose $\mL$ of the form \eqref{L} with $N:=1+p/2$. It satisfies $\mL\in\mathcal L_{AS}^m$ with $m=p/2$, which violates the index condition in Proposition \ref{prop:LASm}, and may hence serve as a counterexample.

When choosing $\mP=\mI$, inequality \eqref{ineq:Lyapunov} is satisfied, and it remains to show that 
\begin{equation}\label{defR}
  \|R(\tau\mL)\|_2 := \sup_{\|u_0\|=1} \|R(\tau\mL) u^0\|\le 1,\quad \mbox{where}\ R(z)=\sum_{j=0}^p \frac{z^j}{j!},
\end{equation}
does \emph{not} hold on any interval $\tau\in[0,\varepsilon)$, with $\varepsilon>0$ arbitrarily small, see Definition \ref{def.stab}(a). Equivalently, we shall show that the matrix function $\mM(\tau):=\mI-R(\tau\mL)^*R(\tau\mL)$ has a negative determinant on $(0,\varepsilon)$ for sufficiently small $\varepsilon>0$. Hence, $R(\tau\mL)^*R(\tau\mL)$ has at least one eigenvalue larger than one, and $\|R(\tau\mL)\|_2>1$ follows.
We shall obtain the inequality $\det \mM(\tau)<0$, for $\tau$ small enough, from the subsequent lemma, thus closing this proof.
$\hfill\square$

\begin{lemma}\label{lem:detM}
Let $p\in4\N$, and $\mL$ of the form \eqref{L} with $N:=1+p/2$. Then, $\mM(\tau):=\mI-R(\tau\mL)^*R(\tau\mL)$ satisfies
\begin{equation}\label{detM-order}
  \det \mM(\tau) =c \tau^{N^2}+\mathcal O(\tau^{N^2+1})\quad \mbox{as } \tau\to0\,,
\end{equation}
with some $c<0$ and $R(z)$ defined in \eqref{defR}.
\end{lemma}
\medskip
\noindent\emph{Proof.}
We only give a sketch of the proof here, the full details are given 
in the Appendix.   

First, we consider the Runge--Kutta method with $p=s=4$, and hence $N=3$. 
In this case, we compute the matrix $\mM(\tau)$ explicitly:
$$
  \mM(\tau) 
    = \left( \begin{array}{rrr}
    \tau^5/12 & \tau^4/4 & \tau^3/3 \\
    \tau^4/4 & 2\tau^3/3 & \tau^2 \\
    \tau^3/3  & \tau^2   & 2\tau
    \end{array}\right) 
    + \left( \begin{array}{rrr}
    \bigO(\tau^6) & \bigO(\tau^5) & \bigO(\tau^4) \\
    \bigO(\tau^5) & \bigO(\tau^4) & \bigO(\tau^3) \\
    \bigO(\tau^4) & \bigO(\tau^3) & \bigO(\tau^2) \\
    \end{array}\right)
    \quad \mbox{as } \tau\to 0\,.
$$
The determinant of the first matrix yields the leading order coefficient:
$$
 \det\mM(\tau) 
=\det \left( \begin{array}{rrr}
	\tau^5/12 & \tau^4/4   & \tau^3/3 \\
	\tau^4/4 & 2\tau^3/3   & \tau^2 \\
	\tau^3/3  & \tau^2     & 2\tau
	\end{array}\right) + \bigO(\tau^{10}) 
= -\frac{\tau^9}{216} + \bigO(\tau^{10})
\quad \mbox{as } \tau\to 0.
$$
Consequently, $\det\mM(\tau)<0$ for sufficiently small $\tau>0$,
which disproves \eqref{defR} for $p=4$.
 
For the general case $p\ge 8$, we insert the stability matrix \eqref{1.R},
$$
  \mM(\tau) = \mI - \bigg(\sum_{j=0}^p\frac{\tau^j}{j!}(\mL^*)^j \bigg)\,
  \bigg(\sum_{\ell=0}^p\frac{\tau^\ell}{\ell!}\mL^\ell \bigg),
$$
expand the matrix coefficients $M_{ij}$ of $\mM(\tau)$ in powers of $\tau$ and use properties of the matrices $\mL$ and $\mL_H$, leading to 
$$
  M_{ij} = \bar{m}_{ij}\tau^{p+3-(i+j)} - \bar{n}_{ij}\tau^{p+1}
  + \bigO(\tau^{p+4-(i+j)}),
$$
where $\bar{m}_{ij}$ and $\bar{n}_{ij}$ are numbers depending on $p$ and $\mL_H$. The determinant of $\mM(\tau)$ can be computed by using the Leibniz formula for determinants and definition $N=p/2+1$:
$$
  \det\mM(\tau) = \det(\bar{m}_{ij}-\bar{n}_{ij})_{1\le i,j\le N}
  \tau^{N^2} + \bigO(\tau^{N^2+1}).
$$
The remaining determinant can be calculated by taking into account the explicit formula of the Hankel determinant:
$$
  c := \det(\bar{m}_{ij}-\bar{n}_{ij})_{1\le i,j\le N}
  = 2^N\bigg(\prod_{i=1}^N\frac{1}{(p/2-i+1)!}\bigg)^2
  \frac{\prod_{1\leq i<j\leq N} (i-j)^2}{\prod_{i,j=1}^N(i+(j-1))}
  \bigg(1 - \binom{p}{p/2}\bigg).
$$
We deduce from $\binom{p}{p/2}>1$ that $c<0$ as claimed in Lemma \ref{lem:detM}. In particular, for sufficiently small $\tau>0$, $\det\mM(\tau)$ is negative, which also finishes the proof of Theorem \ref{th:p4N}.
$\hfill\square$

\medskip
Finally, we note that $\|R(\tau\mL)\|_2>1$, the condition to violate strong stability,  can also be tested numerically. But only for $p=4$ the computation can be carried out in the standard double precision of Matlab: It shows that $\varepsilon=0.304$ can be used, yielding $\max_{[0,\varepsilon]}\|R(\tau\mL)\|_2-1\approx 1.3$E-6. 
For larger values of $p$, this computation is numerically so sensitive that we used octuple precision: For $p=8$ one can use $\varepsilon=0.027$, yielding $\max_{[0,\varepsilon]}\|R(\tau\mL)\|_2-1\approx 9.0$E-22, 
and for $p=12$,  $\varepsilon=0.027$ with $\max_{[0,\varepsilon]}\|R(\tau\mL)\|_2-1\approx 7.2$E-46. 


\section{ACKNOWLEDGMENTS}
The authors were supported by the Austrian Science Fund (FWF) via the FWF-funded SFB \# F65.
The third author acknowledges support from the FWF, grant P33010. 
This work received funding from the European Research Council (ERC) under the European Union's Horizon 2020 research and innovation programme, ERC Advanced Grant NEUROMORPH, no.~101018153.


\nocite{*}
\bibliographystyle{aipnum-cp}%

\section{Appendix: Proof of Lemma \ref{lem:detM}}

Let $p\in4\N$ and let $e_j$ denote the $j$-th Euclidean basis vector. By exploiting the special staircase form of the matrix $\mL$ in \eqref{L}, we find that
$$
 \mL e_1 = e_2, \quad
 \mL e_j = -e_{j-1}+e_{j+1}\quad\mbox{for }j=2,\ldots,N-1, \quad
 \mL e_N = -e_{N-1}-e_N.
$$
A computation shows that, for $j=1,\ldots,N-1$ (with $m=N-1$),
\begin{equation}\label{kerL}
  e_j\in\ker \mL_H, \quad
	\sqrt{-\mL_H}e_j = \sqrt{-\mL_H}\mL e_j = \cdots = \sqrt{-\mL_H}\mL^{m-j}e_j = 0,
	\quad\sqrt{-\mL_H}\mL^{m-j+1}e_j\neq 0
\end{equation}
and $e_N\not\in\ker \mL_H$. 
The coefficients of the matrix $\mM:=I-R(\tau \mL)^*R(\tau \mL)$ are given by 
$M_{ij}=\langle e_i,\mM e_j\rangle$ for $i,j=1,\ldots,N$.
By definition \eqref{1.R} of the stability function and the condition $s=p$, we have
\begin{align*}
 \mM
&= \mI - \bigg(\sum_{j=0}^p\frac{\tau^j}{j!}(\mL^*)^j \bigg)\,
 \bigg(\sum_{\ell=0}^p\frac{\tau^\ell}{\ell!}\mL^\ell 
\bigg) 
\\
&= -\sum_{\ell=1}^{p+1}\frac{\tau^\ell}{\ell!}\sum_{k=0}^\ell\binom{\ell}{k} (\mL^*)^k \mL^{\ell-k} 
+ \frac{\tau^{p+1}}{(p+1)!}\big((\mL^*)^{p+1}+\mL^{p+1}\big) + \bigO(\tau^{p+2})
\quad\mbox{as }\tau\to 0.
\end{align*}
The first term on the right-hand side can be rewritten by using
\cite[(A.6)]{AAC22}: as
\begin{align*}
 \mM 
&= -2\sum_{\ell=1}^{p+1}\frac{\tau^\ell}{\ell!}\sum_{k=0}^{\ell-1} \binom{\ell-1}{k} (\mL^*)^k \mL_H \mL^{\ell-1-k} 
+ \frac{\tau^{p+1}}{(p+1)!}\big((\mL^*)^{p+1}+\mL^{p+1}\big) + \bigO(\tau^{p+2}).
\end{align*}
Therefore, for $i,j=1,\ldots,N$,
\begin{align*}
 M_{ij} 
&= 2\sum_{\ell=1}^{p+1}\frac{\tau^\ell}{\ell!}\sum_{k=0}^{\ell-1} \binom{\ell-1}{k}
\Big\langle\sqrt{-\mL_H} \mL^k e_i,\sqrt{-\mL_H} \mL^{\ell-1-k}e_j\Big\rangle 
\\
&\phantom{xx}{}+ \frac{\tau^{p+1}}{(p+1)!} \big\langle e_i,\big((\mL^*)^{p+1}+\mL^{p+1}\big)e_j\big\rangle + \bigO(\tau^{p+2}) 
=: I_1 + I_2 + \bigO(\tau^{p+2}).
\end{align*}
For $I_1$, due to~\eqref{kerL}, the terms in the sum with $k\leq m-i$ and $\ell-1-k\leq m-j$ vanish, but those with $k\ge m-i+1$ and $\ell-1-k\ge m-j+1$ may not vanish (recall that $m=p/2$). 
Since the $\tau$-order of the leading nonvanishing term equals $\ell=m+2+k-j=2m+3-(i+j)$, we find that
\begin{align*}
 I_1 
&= -\frac{2\tau^{2m+3-(i+j)}}{(2m+3-(i+j))!} \binom{2m+2-(i+j)}{m-i+1} 
\big\langle e_i,(\mL^*)^{m-i+1}\mL_H \mL^{m-j+1}e_j\big\rangle 
+ \bigO(\tau^{2m+4-(i+j)}).
\end{align*}
Because of $p+1=2m+1$, the term $I_2$ is of leading order 
$\tau^{2m+3-(i+j)}$ if and only if $i=j=1$. 
Therefore,
\begin{align*}
 M_{ij} 
&= -\frac{2\tau^{2m+3-(i+j)}}{(2m+3-(i+j))!} \binom{2m+2-(i+j)}{m-i+1}\big\langle e_i,(\mL^*)^{m-i+1}\mL_H\mL^{m-j+1}e_j\big\rangle 
\\
&\phantom{xx}{}+ \delta_{1i}\delta_{j1}\frac{\tau^{p+1}}{(p+1)!} \big\langle e_1,\big((\mL^*)^{p+1}+\mL^{p+1}\big)e_1\big\rangle +\bigO(\tau^{2m+4-(i+j)}),
\end{align*}
where $\delta_{i1}$ is the Kronecker delta.

We claim that the inner products in the first term of $M_{ij}$ are all equal to $-1$. 
The inclusion $\ker \sqrt{-\mL_H} \supseteq \spn\{e_1,\ldots,e_{N-1}\}$ implies that
\begin{align} \label{I1_ij}
 \big\langle e_i&,(\mL^*)^{m-i+1}\mL_H\mL^{m-j+1}e_j\big\rangle
 = -\Big\langle\sqrt{-\mL_H}\mL^{m-i+1}e_i,
 \sqrt{-\mL_H}\mL^{m-j+1}e_j\Big\rangle \nonumber \\
 &= -\Big\langle\sqrt{-\mL_H}\mL^{N-i}e_i,
 \sqrt{-\mL_H}\mL^{N-j}e_j\Big\rangle
 = -\Big\langle \sqrt{-\mL_H}e_N,\sqrt{-\mL_H}e_N\Big\rangle 
 = \langle e_N,\mL_H e_N\rangle= -1\,.
\end{align}
For the second term of $M_{ij}$ we compute similarly: 
\begin{equation}\label{Lp}
 \langle e_1,\mL^{p+1}e_1\rangle
= \langle e_1,\mL^{2m+1}e_1\rangle 
= \langle(\mL^*)^m e_1,\mL\mL^m e_1\rangle
= \langle(\mL^*)^m e_1,\mL_H\mL^m e_1\rangle
= -\langle\sqrt{-\mL_H}(\mL^*)^m e_1,\sqrt{-\mL_H}\mL^m e_1\rangle\,,
\end{equation}
where we used $\langle y,\mL_Sy\rangle=0$ for all $y\in\C^n$. 
We recall from \eqref{I1_ij} that $\sqrt{-\mL_H}\mL^m e_1=e_N$, and obtain analogously that 
$$
  \sqrt{-\mL_H}(\mL^*)^m e_1=(-1)^m e_N=e_N\,,
$$
where we used the property $m=p/2\in 2\N$. Thus \eqref{Lp} equals $-1$. 
So we conclude
\begin{equation*} 
\begin{split}
 M_{ij} 
&= \frac{2\tau^{2m+3-(i+j)}}{(2m+3-(i+j))!} \binom{2m+2-(i+j)}{m-i+1}
- 2\delta_{1i}\delta_{j1}\frac{\tau^{2m+1}}{(2m+1)!}  + \bigO(\tau^{2m+4-(i+j)}) 
\\
&= \bigg( m_{ij}\tau^{2m+3-(i+j)} - 2\delta_{1i}\delta_{j1}\frac{\tau^{2m+1}}{(2m+1)!} \bigg)  +\bigO(\tau^{2m+4-(i+j)}) \,,
 \end{split}
\end{equation*}
where we defined
\begin{equation}\label{mij}
 m_{ij}
 := \frac{2}{(2m+3-(i+j))!}\binom{2m+2-(i+j)}{m-i+1} 
 = \frac{2}{(2m+3-(i+j))(m-i+1)!(m-j+1)!}\,. 
\end{equation}

We proceed with the computation of the determinant of $\mM$.
We use $M_{ij}=\bigO(\tau^{2m+3-(i+j)})$ and Leibniz' formula for determinants:
\begin{equation*}
 \det \mM 
= \sum_{\sigma\in S_N}\sgn(\sigma)\prod_{i=1}^N M_{i,\sigma(i)}
= \sum_{\sigma\in S_N}\sgn(\sigma)\, \bigO(\tau^{\sum_{i=1}^N(2m+3-(i+\sigma(i)))})
= \bigO(\tau^{(m+1)^2}),
\end{equation*}
where $\sgn$ is the sign function of permutations in the permutation group $S_N$ and
\begin{align*}
 \sum_{i=1}^N(2m+3-(i+\sigma(i))) 
 = (m+1)(2m+3) - 2\sum_{i=1}^{m+1} i 
 = (m+1)(2m+3) - (m+1)(m+2) = (m+1)^2.
\end{align*}
This shows that
\begin{equation}\label{detM}
 \det \mM 
= \tau^{(m+1)^2}\det\bigg(m_{ij}-\delta_{1i}\delta_{1j} \frac{2}{(2m+1)!}\bigg)_{1\le i,j\le N}  +\bigO(\tau^{(m+1)^2+1})\,,
\end{equation}
matching the claim \eqref{detM-order}. So it remains to determine the sign of the leading term.

The Laplace expansion with respect to the first column yields for the
determinant on the right-hand side:
\begin{equation}\label{det1}
 \det\bigg(m_{ij}-\delta_{1i}\delta_{1j} \frac{2}{(2m+1)!}\bigg)_{1\le i,j\le N} 
= \det(m_{ij})_{1\le i,j\le N} - \frac{2}{(2m+1)!} \det(m_{ij})_{2\le i,j\le N}\,. 
\end{equation}
The first determinant on the right-hand side becomes, by Leibniz' formula for determinants and definition \eqref{mij} of $m_{ij}$,
\begin{align*}
 \det(m_{ij})_{1\le i,j\le N}
&= \sum_{\sigma\in S_N}\sgn(\sigma)\prod_{i=1}^N \frac{2}{(2m+3-(i+\sigma(i)))(m-i+1)!(m-\sigma(i)+1)!} 
\\ 
&= 2^N\bigg(\prod_{i=1}^N\frac{1}{(m-i+1)!}\bigg)^2\sum_{\sigma\in S_N}
 \sgn(\sigma) \prod_{i=1}^N\frac{1}{2m+3-(i+\sigma(i))} 
\\ 
&= 2^N\bigg(\prod_{i=1}^N\frac{1}{(m-i+1)!}\bigg)^2 \det\bigg(\frac{1}{2m+3-(i+j)}\bigg)_{1\le i,j\le N}. 
\end{align*}
The determinant of the Hankel matrix can be  determined explicitly from Cauchy's double alternant (see, e.g., \cite[Theorem 12]{Kra05}):
\begin{equation*}
  \Delta_1 
:= \det\bigg(\frac{1}{2m+3-(i+j)}\bigg)_{1\le i,j\le N}
 = \det\bigg(\frac{1}{i+(j-1)}\bigg)_{1\le i,j\le N} 
 = \frac{\prod_{1\leq i<j\leq N} (i-j)^2}{\prod_{i,j=1}^N(i+(j-1))}>0\,.
\end{equation*}

The second determinant on the right-hand side of \eqref{det1} is treated similarly:
\begin{equation*}
 \det(m_{ij})_{2\le i,j\le N}
= 2^{N-1}\bigg(\prod_{i=2}^{N}\frac{1}{(m-i+1)!}\bigg)^2\Delta_2, 
\end{equation*}
where
\begin{equation*}
 \Delta_2 
:= \det\bigg(\frac{1}{2m+3-(i+j)}\bigg)_{2\le i,j\le N}
 = \det\bigg(\frac{1}{i+(j-1)}\bigg)_{1\le i,j\le N-1} \\
 = \frac{\prod_{1\leq i<j\leq N-1} (i-j)^2}{\prod_{i,j=1}^{N-1}(i+(j-1))}.
\end{equation*}
In fact, we can express $\Delta_1$ in terms of $\Delta_2$:
\begin{equation*}
 \Delta_1 
= \frac{\prod_{1\leq i<j\leq N-1} (i-j)^2}{\prod_{i,j=1}^{N-1}(i+(j-1))}\
 \frac{(2N-1)\prod_{i=1}^{N-1} (i-N)^2}{\prod_{i=1}^N(i+N-1)^2}
= \frac{\Delta_2}{\Delta_3},
\end{equation*}
where, because of $N=m+1$,
\begin{align*}
 \Delta_3 
 &:= \frac{\prod_{i=1}^N(i+N-1)^2}{(2N-1)\prod_{i=1}^{N-1} (i-N)^2}
 = \frac{\big((2N-1)!/(N-1)!\big)^2}{(2N-1)(N-1)!^2} \\
 &= \frac{(2N-1)!}{(N-1)!^2}\frac{(2N-2)!}{(N-1)!^2}
 = \frac{(2m+1)!}{m!^2}\frac{(2m)!}{m!^2} = \frac{(2m+1)!}{m!^2}\binom{2m}{m}.
\end{align*}
Therefore, \eqref{det1} becomes
\begin{align*}
 c:=\det\bigg(m_{ij}-\delta_{1i}\delta_{1j} \frac{2}{(2m+1)!}\bigg)_{1\le i,j\le N} 
&= 2^N\bigg(\prod_{i=1}^N\frac{1}{(m-i+1)!}\bigg)^2\Delta_1 -\frac{2^N}{(2m+1)!}\bigg(\prod_{i=2}^{N}\frac{1}{(m-i+1)!}\bigg)^2 \Delta_2
\\
&= 2^N\bigg(\prod_{i=1}^N\frac{1}{(m-i+1)!}\bigg)^2\Delta_1  \bigg(1 - \frac{(m!)^2}{(2m+1)!}\Delta_3\bigg) 
\\
&= 2^N\bigg(\prod_{i=1}^N\frac{1}{(m-i+1)!}\bigg)^2\Delta_1  \bigg(1 - \binom{2m}{m}\bigg)<0\,.
\end{align*}
Thus, \eqref{detM} shows that, for sufficiently small $\tau>0$, the sign of $\det \mM$ is 
\begin{equation*}
  \sgn(\det \mM) = \sgn\bigg(1 - \binom{2m}{m}\bigg) <0,
\end{equation*}
finishing the proof. $\hfill\square$

\end{document}